\documentclass[a4paper,11pt]{amsart}
\usepackage{a4}
\usepackage{amsfonts,amssymb}
\usepackage{amsmath}
\usepackage{amsopn}
\usepackage{amsthm}
\usepackage{latexsym}

\newtheorem{theorem}{Theorem}
\newtheorem{lemma}[theorem]{Lemma}
\newtheorem{prop}[theorem]{Proposition}
\newtheorem{cor}[theorem]{Corollary}
\newtheorem{rmk}[theorem]{Remark}

\newtheoremstyle{Probleme}% name
  {3pt}%      Space above, empty = `usual value'
  {3pt}%      Space below
  {\it}% Body font
  {}%         Indent amount (empty = no indent, \parindent = para indent)
  { }% Thm head font
  { }%        Punctuation after thm head
  { }%     Space after thm head: " " = normal interword space;
  {}% Thm head spec
\theoremstyle{Probleme}
\newtheorem*{Prob}{{\bf \sf (P)}}

\newcommand {\be }{ \begin{equation} }
\newcommand {\bea }{ \begin{eqnarray} }
\newcommand {\eea }{ \end{eqnarray} }
\newcommand {\ba }{\begin{array}}
\newcommand {\ea} {\end{array}}

\newcommand {\lbl }{\label}

\newcommand {\non }{\nonumber  }

\newcommand {\Z }{\mathbb{Z}}
\newcommand {\ml } {\left( \begin{array}}
\newcommand {\mr} {\end{array} \right) }

\newcommand\id{\operatorname{id}}
\newcommand\Gal{\operatorname{Gal}}
\newcommand\h{\operatorname{\mathcal{H}}}
\newcommand\Tr{\operatorname{Tr}}
\newcommand\can{\operatorname{can}}

\newcommand\Comod{\operatorname{Comod}}

\newcommand {\Oq} {\mathcal{O}_q(SL(2))}
\newcommand {\B} {\mathcal{B}}

\newcommand {\BB} {\Box _{\B(E)}}

\newcommand {\nedots}{ \cdot^{\cdot ^{\cdot} }}

\title[On the classification of Galois objects]{On the classification of Galois objects over the quantum group of a nondegenerate bilinear form}

\date{\today}
\author{Thomas AUBRIOT}
\address{Institut de Recherche Math\'ematique Avanc\'ee \\ 
C.N.R.S. - Universit\'e Louis Pasteur \\ 
7 rue Ren\'e Descartes \\ 
67084 Strasbourg Cedex, France \\ 
Fax : +33 (0)3 90 24 03 28}

\email{aubriot@math.u-strasbg.fr}

\thanks{The author thanks his advisors J. Bichon and C. Kassel for their day-to-day help.}

\begin{document}
\begin{abstract} We study Galois and bi-Galois objects over the quantum group of a nondegenerate bilinear form, including the quantum group $\Oq$. 
We obtain the classification of these objects up to isomorphism and some partial results for their classification up to homotopy. \end{abstract}


\maketitle
\section*{Introduction}Hopf-Galois extensions and objects are quantum analogues of principal fibre bundles and torsors. It is in general a 
difficult problem to classify these objects. Several authors have already contributed to this problem, mainly in the finite dimensional 
case (\cite{Ma1}, \cite{Ma2}, \cite{Sa3}, \cite{PO}, \cite{B2}...). In this paper we study a very different class of infinite dimensional Hopf 
algebras, including the quantum group~$\Oq$ of functions over~$SL_2$. We obtain the classification up to isomorphism and present some partial 
results for the classification up to homotopy. Homotopy for Hopf-Galois extensions was introduced by Kassel~\cite{K1} and developped with 
Schneider \cite{KS} in order to classify Galois extensions up to a coarser equivalence relation than isomorphism. This relation is very useful for 
pointed Hopf algebras but it appears that, when the Hopf algebra is the quantum group~$\Oq$, the classification up to homotopy is harder to obtain 
than the one up to isomorphism.

We consider the Hopf algebras $\B(E)$ introduced by Dubois-Violette and Launer \cite{DV-L} as the quantum groups of nondegenerate 
bilinear forms  given by invertible matrices~$E$ over a field~$k$. One simple and interesting example of such a Hopf algebra is the 
quantum group $\Oq$ of functions over~$SL_2$. Bichon~\cite{B} has proved that the representation category of each Hopf algebra~$\B(E)$ 
is monoidally equivalent to the one of $\Oq$, where~$q$ is a solution of the equation
\begin{equation}q^2 + \Tr(E^{-1}E^t) q +1 =0.\non \end{equation}
The main ingredient of his proof is the construction of a $\B(E)$-$\Oq$-bi-Galois object~$\B(E,E_q)$ for a well-chosen 
invertible matrix~$E_q$. In fact, such Galois objects~$\B(E,F)$ can be defined even when $k$ is only assumed to be a commutative ring. 
They are generic in the following sense: if~$k$ is a PID (principal ideal domain), for any~$\B(E)$-Galois object~$Z$ there exist 
an integer~$m \geq 2$ and an invertible matrix~$F \in GL_m(k)$ such that $Z$ is isomorphic to $\B(E,F)$. Moreover, two 
such~$\B(E)$-Galois objects~$\B(E,F_1)$ and~$\B(E,F_2)$ are isomorphic if and only if there exists~$P\in GL(k)$ such 
that~$F_1 =P F_2 P^t$ ( $P^t$ denotes the transpose of~$P$). In the case when~$k$ is a field, we obtain a full classification up to isomorphism of the Galois objects of~$B(E)$. 
As a consequence, the group of~$\B(E)$-bi-Galois objects is trivial as well as the lazy cohomology group. 
Note that Ostrik~\cite{O} has recently classified module categories over representations of~$\Oq$, which 
together with Ulbrich's and Schauenburg's work (see \cite{U1}, \cite{U2} and \cite{S2}) also yields a classification of Galois objects, but the tools 
used in~\cite{O} are very different from ours. 

Concerning the classification up to homotopy, we prove a partial result. Namely, we show that two Galois objects~$\B(E,F_1)$ and~$\B(E,F_2)$ are
 homotopically equivalent if the matrices~$F_1^{-1} F_1 ^t$ and~$F_2^{-1} F_2^t$ have the same characteristic polynomial. In particular, any cleft~$\Oq$-Galois 
object is homotopically trivial. 

The paper is organized as follows. In Section~$1$ we recall some basic facts on Galois and bi-Galois extensions. Section~$2$ and~$3$ are 
devoted to the isomorphism problem for~$\B(E)$-Galois objects, while Section~$4$ deals with the classification up to homotopy. 

\section{Hopf-Galois extensions and bi-Galois objects}
Let $k$ be a commutative ring. All objects in this paper
belong to the tensor category of~$k$-modules and the tensor product over $k$ is denoted by $\otimes $.
Let~$H$ be a Hopf algebra and~$Z$ be a left $H$-comodule algebra with coaction $\delta : Z \to H \otimes Z $. We define the
subalgebra $R=Z^{co H}$ of {\bf $H$-coinvariant elements} of $Z$ by
\be R = \{ z\in Z \mid \delta(z) = 1\otimes z\}. \non \end{equation} The linear application $\can 
: Z \otimes _R Z \to H \otimes Z $ given by \be \can(z\otimes z')
= \delta(z)(1\otimes z') \non \end{equation} for all $z$, $z'\in Z$, is called the
{\bf canonical map} of~$Z$.

If~$Z$ is a left~$H$-comodule algebra and~$R$ is a subalgebra of~$Z$, then we say that~$R\subset Z$ is a {\bf $H$-Galois extension} if the 
subalgebra of $H$-coinvariant elements 
is~$R$ and if the canonical map~$\can : Z \otimes_R Z \to H \otimes Z$ of~$Z$ is an isomorphism. In this case, we also say that~$Z$ is an~$H$-Galois 
extension of~$R$. 
A Galois extension~$Z$ of~$R$ is said to be {\bf faithfully flat} if~$Z$ is faithfully flat as a right or left $R$-module. 
An {\bf $H$-Galois object} is an $H$-Galois extension of $k$ which is $k$-faithfully flat.

A {\bf morphism of Galois extensions} between two $H$-Galois
extensions~$Z$ and~$Z'$ of~$R$ is a morphism of~$H$-comodule algebras which is the identity on~$R$. If~$Z'$ is faithfully flat, it is 
always an isomorphism. We denote~$\Gal _R(H / k)$ the set of isomorphism classes of faithfully flat~$H$-Galois extensions of~$R$. 
If~$Z$ is a faithfully flat~$H$-Galois extension of~$R$, its isomorphism class in~$\Gal _R(H/k)$ is denoted by~$[Z]$. If one of the objects~$R$ or~$k$ is clear, we will 
omit it from the notation. In the same way, one can define right~$H$-Galois extensions of~$R$ and we denote~$\Gal ^r_R(H / k)$ the set of 
isomorphism classes of faithfully flat right~$H$-Galois extensions. If~$H$ has a bijective antipode, then there is a bijection between the sets~$\Gal _R(H/k)$ 
and~$\Gal _R ^r(H/k)$. 

Recall that, if $H$ is a Hopf algebra, $U$ is a right $H$-comodule and $V$ a left
$H$-comodule, the {\bf cotensor product}~$U \Box _H V$ is the kernel
of the map 
\begin{equation} \delta_U \otimes \id _V
-\id _U \otimes \delta _V : U\otimes V \to U\otimes
H \otimes V , \non \end{equation} (or the equalizer of the coactions of $U$ and
$V$). 
%In particular, since the cotensor product $U\Box _H V $ of $U$
%and $V$ over $H$ is contained in $U \otimes V$, we use the notation~$u\otimes v$ for its elements. 

A bilinear map $\sigma : H\times H \rightarrow k$ is a {\bf right invertible cocycle} for the Hopf algebra $H$ if $\sigma$ is 
convolution-invertible and satisfies the relations  
\be \lbl{left_cocycle} 
 \sigma (x_{(1)}y_{(1)},z) \sigma (x_{(2)},y_{(2)})=\sigma (x,y_{(1)}z_{(1)}) \sigma  (y_{(2)},z_{(2)}) \non 
\end{equation} and \be \lbl{norm} \sigma  (1,x)= \sigma (x,1)=\varepsilon (x), \non \end{equation} for all~$x,y,z \in H$. Here~$\varepsilon$ denotes the 
counit of~$H$ and we have used Sweedler's notation~$x_{(1)}\otimes x_{(2)}$ for the comultiplication. 
Note that we use right cocycles whose definition is different from the one of left cocycles (see~\cite{M}). We denote~$\sigma ^{-1} $ the inverse 
of~$\sigma$ for the convolution ; $\sigma^{-1}$ is a left cocycle. 

Recall (\cite[Chapter 7]{M}) that if $H$ is a Hopf algebra, $\sigma :H\times H \rightarrow k $ an invertible cocycle, one can 
define the Hopf algebra~$H^{\sigma}$ as the coalgebra $H$ with the twisted product 
\be
\lbl{produit_tordu_hopf} x \cdot _{\sigma} y = \sigma ^{-1} (x_{(1)} ,y_{(1)})  x_{(2)} y_{(2)} \sigma  (x_{(3)} ,y_{(3)}) \non \end{equation} 
and the $H$-comodule algebra $H_{\sigma}$ as the left $H$-comodule $H$ with the twisted product 
\be
\lbl{produit_tordu_clive_gauche} x \cdot _{\sigma} y =   x_{(1)} y_{(1)} \sigma (x_{(2)} ,y_{(2)}) , \non \end{equation} for any $x,y \in H$. 
The $H$-comodule algebra $H_{\sigma}$ is an~$H$-Galois extension of~$k$ and all such Galois extensions are called {\bf cleft 
Galois extensions}. If~$H$ is $k$-faithfully flat, it is a {\bf cleft Galois object}.

Kassel and Schneider \cite{KS} (see also \cite{K1}) have defined an
equivalence relation denoted $\sim$ and called {\bf homotopy} on the
class of faithfully flat Galois extensions of~$R$. Two Hopf-Galois extensions are homotopy equivalent if 
there exists a polynomial path between these extensions. More precisely, let $k[t]$ be the algebra of 
polynomials with coefficients in the ground ring $k$. For any~$k$-module~$V$, we denote~$V[t]= V \otimes k[t]$ and 
for~$i\in \{0,1\}$ we denote~$[i]: V[t] \to V $ the~$k$-linear map sending~$vt^n$ to~$vi^n$. These two maps~$[i]$ induce two maps 
\begin{equation} [i]_* : \Gal _{R[t]}(H[t], k[t]) \to\Gal _{R}(H, k) ,\non \end{equation} 
for~$i=0,1$. We say that two~$H$-Galois extensions~$Z_0$ and~$Z_1 \in \Gal _R (H/k)$ are 
{\bf homotopy equivalent} if there exists~$Z \in \Gal _{R[t]} (H[t] /k[t])$ such that~$[i]_* (Z)=Z_i $ for~$i\in \{0,1\}$. 
We denote~$\h _R (H) $ the set of homotopy classes of faithfully flat left~$H$-Galois extensions of~$R$.  

Kassel and Schneider \cite[Proposition 1.6, Corollary 1.11]{KS} have proved that twists of homotopy equivalent 
Galois objects are still homotopy equivalent. In fact, the twist is a particular case of the cotensor product by 
a bi-Galois object. We generalize this result now. 

Let $H$ and $K$ be Hopf algebras. An {\bf $H$-$K$-bi-Galois object} is a $H$-$K$-bi-comodule algebra~$Z$ which is a Galois object 
with respect to the right and the left coactions. By work of Schauenburg \cite{S1} (see also \cite{S2}), the set of bi-Galois objects is a groupoid 
with the multiplication given by the cotensor product. In particular, when $H=K$, the cotensor product over $H=K$ puts a structure of 
group on the set of isomorphism classes of~$H$-$H$-bi-Galois objects. 
If~$Z$ is an~$H$-$K$-bi-Galois object, the cotensor product yields a bijective 
map $\varphi _Z: \Gal _k(K) \to \Gal _k(H)$ defined by 
\begin{equation}
\varphi _Z([A])= [Z\Box _K A] \non \end{equation}
for any left~$K$-Galois object~$A$ (see \cite{S1} and \cite{S2} for details). 

\begin{prop}\lbl{prop_map}For any~$H$-$K$-bi-Galois object~$Z$, the map~$\varphi _Z$ induces a bijective 
map~$\overline{\varphi _Z} : \h _k(K) \to \h _k(H)$ between the homotopy classes of left~$K$-Galois objects and of left~$H$-Galois objects. 

%Furthermore, this map~$\varphi$ induces a map~$\bar{\varphi}: \Cleft (H) \to \Cleft (K)$ between the cleft~$H$-Galois objects and the cleft~$K$-Galois 
%objects if and only there exists a cocycle~$\sigma : H\otimes H \to k$ such that~$K=H^{\sigma}$. In that case, the bi-Galois object~$Z$ is isomorphic to~
%$H_{\sigma}$.  
\end{prop}
\begin{proof}
Let~$A_0, A_1$ be homotopically equivalent~$H$-Galois objects via the~$H[t]$-Galois object~$A$. Then~$Z[t]=Z\otimes k[t]$ is an~$H[t]$-$K[t]$-bi-Galois 
object and $Z[t] \Box _{K[t]} A$ is an homotopy between~$Z \Box _K A_0$ and~$Z\Box _K A_1$. There exists a~$K$-$H$-bi-Galois object~$Z^{-1}$ 
inverse of~$Z$ for the groupoid structure of bi-Galois objects and the map~$\overline{\varphi _{Z^{-1}}} : 
\h (H) \to \h (K)$ induced by~$Z^{-1}$ is the inverse of the map~$\overline{\varphi _{Z}} : \h (K) \to \h (H)$ induced by~$Z$. 
\end{proof}

\section{The Hopf algebra $\B (E)$ and the comodule algebra $\B(E,F)$}
Let $k$ be a commutative ring, $n \geq1$ an integer and $E =(E_{ij}) _{1\leq i,j\leq n} \in GL_n(k)$. Following \cite{DV-L}, we define~$\B _k(E)$ (or~$\B (E)$ when the base ring is clear) as the~$k$-algebra generated by~$n^2$ variables~$a_{ij},\; 1\leq i,j\leq n $, submitted to the matrix relations 

\begin{equation} \lbl{def_BE} E^{-1}a^t E a =\mathrm{I_n} = aE^{-1}a^t E , \non \end{equation}
where~$E^{-1}$ is the inverse matrix of~$E$, $a$ is the matrix $(a_{ij})$, $\mathrm{I_n}$ the identity matrix of size $n$ and 
$ a^t$ denotes the transpose of the matrix~$a$. 

The algebra~$\B (E)$ is a Hopf algebra with comultiplication $\Delta$ defined by 
\begin{equation} \Delta (a_{ij}) =\sum _{k=1}^n a_{ik}\otimes a_{kj}, \non \end{equation}  
counit $\varepsilon $ defined by $\varepsilon (a_{ij})=\delta _{ij}$, for any~$i,j=1 ,\ldots ,n$, where~$ \delta _{ij}$ is Kronecker's symbol, 
and antipode~$S$ defined by the matrix identity~$S(a)= E^{-1 }a^t E$.

Note that if~$n=1$, the Hopf algebra~$\B(E)$ is isomorphic to~$k[\Z / 2\Z]$, whose Galois objects are~$k[\Z / 2\Z]_{\sigma}$ 
for~$\sigma \in  H^2 (\Z /2\Z ,k^*) $. In the following, we will only consider the cases where~$n\geq 2$. In this case, this Hopf algebra is the quantum group of the bilinear (but non necessarily symmetric) form defined by the matrix $E$, in the sense that~$\B(E)$ is the universal Hopf algebra such that the bilinear form is a comodule map (for details see \cite{DV-L}).
% in the following sense (for details see \cite{DV-L}). 
%The {\bf fundamental comodule of~$\B(E)$}, denoted $V_E$, is the finite free~$k$-module of rank~$n$ with its canonical basis~$(v_i)$ and endowed 
%with the~$\B(E)$-comodule structure defined by~$\delta (v_i) = \sum _{k=1}^n v_k \otimes a_{ki}$ for~$1\leq i\leq n$. Then the linear map~$\beta _E : 
%V_E\otimes V_E \to k$ defined by~$\beta _E (v_i ,v_j)= E_{ij}$ for~$1\leq i,j\leq n$ is a~$\B(E)$-comodule morphism and the Hopf algebra is universal for %this property: let~$H$ be a Hopf algebra and ~$V$ an~$H$-comodule that is a free $k$-module of rank~$n$. Let~$\beta : V \otimes V \to k$ be an~$H$-
%comodule morphism such that the associate bilinear form is non-degenerate. Then there exists~$E\in GL_n(k)$ such that~$V$ is a~$\B(E)$-comodule, 
%such that~$\beta$ is a~$\B(E)$-comodule morphism and there exists an unique Hopf algebra morphism~$\varphi : \B(E) \to H$ such that~$(\mathrm{Id}
% \otimes  \varphi ) \circ \delta = \delta '$, where~$\delta $ and~$\delta '$ denote the coactions on~$V$ of~$\B(E)$ and~$H$ respectively.

If $q \in k^*$ is an invertible element of the ring~$k$, let $E_q\in GL_2(k)$ be the matrix defined by 
\be \lbl{def_Eq} E_q =\ml{cc} 0&1\\-q^{-1} &0 \\ \mr . \non 
\end{equation} 
The Hopf algebra $\B(E_q)$ is isomorphic to the Hopf algebra $\Oq$ (see the definition of $\Oq$ in \cite{K2}). 

Let $n,m\geq1 $ be integers and let $E\in GL_n(k)$, $F\in GL_m(k)$ be invertible scalar matrices. Following Bichon \cite{B}, we define the 
algebra $\B(E,F)$ as the $k$-algebra generated by~$n\times m$ variables~$z_{ij}, i=1,\ldots ,n; j=1,\ldots ,m$, submitted to the matrix relations 
\begin{equation} \lbl{def_BEF} F^{-1}z^t Ez =\mathrm{I_m}, \quad zF^{-1}z^tE =\mathrm{I_n} , \non \end{equation} 
where $z$ is the matrix of generators $z_{ij}$ and $\mathrm{I_m,I_n}$ are the identity matrices of size~$m,n$ respectively. 
We consider the $k$-algebra morphism 
$\delta :\B(E,F)\to \B(E) \otimes \B (E,F) $, defined by 
\begin{equation} \lbl{def_coaction_BEF} \delta (z_{ij}) = \sum _{k=1}^n a_{ik}\otimes z_{kj} , \end{equation}
for any $i=1,\ldots ,n$ and $j=1,\ldots ,m$, that endows~$\B(E,F)$ with a left~$\B(E)$-comodule algebra structure. 

In the same way, we have a~$k$-algebra map~$\rho : \B(E,F) \to \B(E,F) \otimes \B(F)$ defined by 
\begin{equation}  \rho (z_{ij}) = \sum _{k=1}^n z_{ik}\otimes b_{kj} ,\non  \end{equation} where the~$b_{ij}$'s stands for the canonical generators of~$\B(F)$. 
The algebra morphism~$\rho$ endows~$\B(E,F)$ with a right comodule algebra structure and~$\B(E,F)$ is a~$\B(E)$-$B(F)$-bicomodule algebra. 

%Similarly to \cite{B}, we endow~$\B(E,F)$ with a $\B(E)$-$\B(F)$-bi-comodule algebra structure. 
Bichon has proved \cite[Propositions 3.3, 3.4]{B} that if~$k$ is a field and if $\Tr(E^{-1}E^t)= \Tr(F^{-1}F^t)$, then $\B(E,F)$ is a $\B(E)$-$\B(F)$-bi-Galois object. Note that the matrices of form $F^{-1}F^t$ appear in Riehm's work \cite{R} on the classification of bilinear form. Precisely, for any nondegenerate bilinear map~$\beta : V \times V \to k$ given by an invertible matrix~$F$, the matrix~$\sigma =F^{-1}F^t$ is called the {\bf asymmetry} of~$\beta$. Over a commutative ring, Bichon's result extends to the following proposition.

\begin{prop}\lbl{prop_galois_object}
The canonical map of~$\B(E,F)$ considered as a left (resp right) $\B(E)$-comodule algebra (resp~$\B(F)$-comodule algebra) is bijective. 
%If~$k$ is a commutative ring, the $k$-algebra $\B  (E,F)$ has a left~$\B (E)$-comodule algebra structure with bijective canonical map $\can : \B(E,F) 
%\otimes  \B(E,F) \to \B(E) 
%\otimes \B(E,F)$ defined as usual by 
%\begin{equation} \lbl{def_can} \can (x\otimes y)= \delta (x) (1 \otimes y),\end{equation} for any $x,y \in  \B (E,F)$. In the same way, the algebra~$\B(E,F)$ 
%has a compatible 
%structure of right~$\B(F)$-comodule algebra with bijective canonical map. 

Moreover, if~$\B (E,F)$ is~$k$-faithfully flat, it is a~$\B(E)$-$\B(F)$-bi-Galois object. 
\end{prop}
\begin{proof}
The proof is the same as for \cite[Propositions 3.3, 3.4]{B}. \end{proof}

Together with Proposition \ref{prop_map}, this yields the following corollary. 

\begin{cor} Assume that~$k$ is a field. 
Let $E$ be an invertible matrix and $q\in k^*$ such that $\mathrm{Tr}(E^{-1}E^t)= -q-q^{-1}$, then there is a bijection between~$\h _k(\B(E))$ and~$\h _k(\Oq )$. 
\end{cor}

\section{Classification up to isomorphism} 
The ~$\B(E)$-comodule algebras~$\B(E,F)$ are generic in the following sense. 

\begin{theorem}\lbl{theorem_classif_iso}
Let $k$ be a PID, $n\geq 2$ an integer, $ E \in GL_n(k)$ and~$Z$ be a~$\B(E)$-Galois object. 
Then there exist an integer~$m\geq 2$ and an invertible matrix~$F\in GL_m(k)$ such that~$\Tr(F^{-1}F^t) = \Tr(E^{-1}E^t)$ 
and such that~$Z$ is isomorphic to~$\B(E,F)$ as a~$\B(E)$-Galois object. 
\end{theorem} 

\begin{proof}Let $\Comod$-$\B(E)$ be the monoidal category of right~$\B(E)$-comodules, with the tensor product $\otimes $ over $k$, 
and $\mathrm{Mod} (k)$ be the monoidal category of~$k$-modules. Following Ulbrich \cite{U1}, \cite{U2} and Schauenburg~\cite{S2}, to 
any~$\B(E)$-Galois 
object $Z$, we associate the fibre functor~$\omega _Z : \Comod $-$\B(E) \to \mathrm{Mod} (k)$ defined by 
\begin{equation} \lbl{def_foncteur_fibre}\omega _Z (V) = V \Box _{\B(E)} Z \non 
 \end{equation} 
for any $V \in \Comod$-$\B(E) $. The map~$Z \to \omega _Z$ defines a bijective correspondence between the left~$\B(E)$-Galois objects 
(they are by definition faithfully flat) and the exact monoidal functors (= fibre functors) $\Comod$-$B(E) \to \mathrm{Mod}(k)$. 
Moreover, the fibre functor $\omega _Z$ sends comodules that are finitely generated projective~$k$-modules to finitely generated projective~$k$-modules (the functor~$\omega _Z$ preserves the duals). We denote by $\psi _2 : \omega _Z (V) \otimes \omega _Z (V) \to \omega _Z (V\otimes V)$ and $\psi_0 : \omega _Z (k) \to k$ the monoidal isomorphisms. Note that~$\psi _2 : (V\Box _{\B(E)} Z )\otimes (V\Box _{\B(E)} Z) \to  (V\otimes V) \Box _{\B(E)} Z$ is induced by the multiplication of~$Z$.

The {\bf fundamental comodule of~$\B(E)$}, denoted $V_E$, is the finite free~$k$-module of rank~$n$ with basis~$(v_1, \ldots ,v_n)$ and endowed 
with the~$\B(E)$-comodule structure defined by~$\delta (v_i) = \sum _{k=1}^n v_k \otimes a_{ki}$ for~$1\leq i\leq n$. 
The linear map~$\beta _E : V_E\otimes V_E \to k$ defined by~$\beta _E (v_i ,v_j)= E_{ij}$ for~$1\leq i,j\leq n$ is a~$\B(E)$-comodule morphism and 
induces a map
\begin{equation} \overline{\beta _E}: W \otimes W \xrightarrow{\psi _2 } \omega _Z (V_E \otimes V_E ) \xrightarrow{\omega _Z (\beta _E)} \omega _Z (k) \xrightarrow{\psi _0} k, \non \end{equation} where $W = \omega _Z(V_E)$. Since~$V_E$ is free of finite rank, $W$ is a finitely generated projective $k$-module. The base ring~$k$ being principal, it implies that~$W$ is a free $k$-module of finite rank, say~$m$. 

Set~$F_{ij}= \overline{\beta _E }(w_i \otimes  w_j)$ for all $1\leq i,j\leq m$. Writing the elements~$(w_{j})_{1\leq j\leq m}$ as elements of~$V_E \otimes Z$ and expanding them in the basis~$(v_1,\ldots ,v_n)$ of~$V _E$, we see that there exist~$(t_{ij})_{i=1,\ldots ,n; j=1,\ldots ,m} \in Z$ such that~$w_j= \sum _{i=1}^n v_i  \otimes t_{ij}$ for any $j=1,\ldots ,m$. Since~$(w_j)_{1\leq i\leq m}$ belong to the cotensor product~$V_E \Box_{\B(E_q)} Z$, the elements~$(t_{ij})_{i=1,\ldots ,n; j=1, \ldots ,m}$ satisfy the relations 
\be \lbl{coaction_T}\delta (t_{ij})= \sum _{k=1}^n a_{ik} \otimes t_{kj}\end{equation} for all $1\leq i \leq n$ and $1\leq j\leq m$. 

Since the monoidal isomorphism~$\psi _2$ is given by the multiplication of~$Z$, the image of the base $(w_j)_{1\leq j\leq m}$ by the map $\overline{\beta_E}$ is equal to
\begin{equation} \begin{array}{rcl}
F_{ij} &=& \overline{\beta _E}(w_i \otimes   w_j) \\
&= & \psi _0  \circ (\beta _E \otimes \id) \circ \psi _2 ((\sum _{k=1}^n v_k \otimes  t_{ki}) \otimes  (\sum _{l=1}^n v_l \otimes  t_{lj})) \\
&= & \psi _0  \circ (\beta _E\otimes \id) (\sum _{k,l=1}^n ( v_k  \otimes   v_l) \otimes t_{ki} t_{lj}) \\
&= & \psi _0 (\sum _{k,l=1}^n E_{kl} \otimes t_{ki} t_{lj}) \\
&=& \sum _{k,l=1}^n E_{kl}  t_{ki} t_{lj} \\  \end{array} \non  \end{equation}
for any $1\leq i,j\leq m$. Putting $T= (t_{ij})_{i=1,\ldots ,n;j=1,\ldots ,m}$ and $F=(F_{ij})_{1\leq i,j\leq m}$, we obtain
\be \lbl{rel_bar1} F=T^t ET.\end{equation} 

Let us now consider the $k$-linear map $\nu : k \to V_E \otimes V_E $ defined by 
\begin{equation} \lbl{def_nu} \nu (1) = \sum _{i,j =1,\ldots ,n} E^{-1}_{ij} v_i \otimes v_j , \non 
\end{equation} 
where $E^{-1}_{ij}$ denotes the $(i,j)$-entry of the inverse matrix $E^{-1}$. Since this map is a~$\B(E)$-comodule morphism, it induces a linear map
\begin{equation} \bar{\nu}: k \xrightarrow{\psi _0 ^{-1} } \omega _Z (k) \xrightarrow{\omega _Z (\nu)} \omega _Z(V_E \otimes V_E) \xrightarrow{\psi _2 ^{-1}} \omega _Z (V_E) \otimes \omega _Z (V_E). \non \end{equation} 
Let us compute $\bar{\nu}(1)$. We have

\begin{equation} \begin{array}{rcl}\bar{\nu} (1) &=& \psi _2^{-1} \circ (\nu \otimes \id ) \circ \psi _0 ^{-1} (1) \\
&=&  \psi _2^{-1} \circ (\nu \otimes \id ) (1\otimes 1) \\
&=& \psi _2^{-1} (\sum _{i,j=1}^n E^{-1}_{ij} (v_i \otimes v_j ) \otimes 1 \\
&=& \sum _{k,l=1}^n E^{-1}_{kl} (v_k \otimes 1 ) \otimes (v_l  \otimes 1) . \\ \end{array} \non \end{equation} 
 
Expanding $\bar{\nu}(1)$ in the basis $(w_j)_{1\leq j\leq m}$ of $\omega _Z(V_E)$, we obtain a matrix~$(G_{ij})_{1\leq i,j \leq m}\in M _m(k)$ such that \be \bar{\nu}(1)=\sum _{i,j=1}^m G_{ij} w_i \otimes  w_j =\sum _{i,j=1}^m \sum _{k,l=1}^n G_{ij}(v_k \otimes t_{ki}) \otimes  (v_l \otimes  t_{lj}) \non  \end{equation} 
for all $1\leq i,j\leq m$. Then we have 
\begin{equation} 
 \sum _{k,l=1}^n E^{-1}_{kl} (v_k \otimes 1 ) \otimes (v_l  \otimes 1) = \sum _{i,j=1}^m \sum _{k,l=1}^n  G_{ij} (v_k \otimes t_{k,i}) \otimes  (v_l \otimes t_{lj}) ,\non  \end{equation} and then 
\begin{equation} E^{-1}_{kl} = \sum _{i,j=1}^m  G_{ij} t_{ki} t_{lj} ,\non \end{equation} 
which we can rewrite as
\begin{equation} \lbl{rel_bar2}  E^{-1}=TGT^t .\end{equation} 

We now prove $G=F^{-1}$. We have 
\begin{equation} (\beta _E \otimes \id _{V_E}) \circ (\id _{V_E} \otimes \nu ) = \id _{V_E}\non  \end{equation} and 
\begin{equation} (\id _{V_E} \otimes \beta _E ) \circ (\nu \otimes \id _{V_E}) = \id _{V_E}\non   \end{equation} 
for $\beta _E $ and $\nu$. Since $\omega _Z$ is monoidal, we obtain 
\begin{equation} (\overline{\beta _E} \otimes  \id _W) \circ (\id _W \otimes  \bar{\nu} ) = \id _W\non \end{equation} and 
\begin{equation} (\id _W \otimes  \overline{\beta _E} ) \circ (\bar{\nu} \otimes  \id _W) =\id _W \non \end{equation} for $\overline{\beta _E}$ and $\bar{\nu}$. 
That is for any basis vector $w_i$ we have 
\begin{equation} \sum _{jk}^m F_{ij} G _{jk} w_k =w_i \quad \mathrm{and}\quad  \sum _{jk}^m G_{jk} F_{ki} w_j =w_i.\non  \end{equation} 
This implies that the matrix $G$ is the inverse of $F$. Then Relations~(\ref{rel_bar1}) and~(\ref{rel_bar2}) yield the relations 
\begin{equation} \lbl{rel_FT}  F^{-1}T^tET=\mathrm{I}_m \quad \mathrm{and} \quad TF^{-1}T^tE=\mathrm{I}_n .\end{equation}

In the same way, we obtain
\begin{equation} \beta _E \circ \nu (1) = \mathrm{Tr}(E^{-1}E^t ) .\non \end{equation} 
Since $\omega _Z$ is monoidal, 
\begin{equation} \overline{\beta _E} \circ \bar {\nu} (1) =  \Tr(F^{-1}F^t )\non \end{equation} 
has to be equal to $ \Tr(E^{-1}E^t ) $. 
When~$\bar{k}$ is a field, Bichon has proved in~\cite[Section 4]{B} that, under this condition, the algebra~$\B _{\bar{k}} (E,F)$ is nonzero. Since our base ring~$k$ is a PID, it embeds into a field~$\bar{k}$. It is clear that for any invertible matrices $E,F$, the algebras~$\B _k (E,F)\otimes _k \bar{k}$ and~$\B _{\bar{k}} (E,F)$ are isomorphic. Therefore, $\B _k(E,F)$ is nonzero provided~$\Tr (E^{-1}E^t )=\Tr(F^{-1}F^t )$. 

In view of~(\ref{rel_FT}) the map
\be \varphi (z_{ij})=t_{ij},\non  \end{equation} 
defines an algebra morphism~$\varphi : \B(E,F) \to Z$. We claim that~$\varphi $ is an isomorphism of~$\B(E)$-Galois objects. 
First to see that~$\varphi $ is a $\B(E)$-comodule morphism, it is enough to check it on the generators~$(z_{ij})$. The definition of the coaction~(\ref{def_coaction_BEF}) and relation~(\ref{coaction_T}) give
 \begin{equation} ( \mathrm{Id}\otimes \varphi )\circ \delta _{\B (E,F)} (z_{ij}) = \sum _{k=1} ^n a_{ik}\otimes t_{kj}= \delta _Z \circ \varphi (z_{ij})\non  \end{equation} for any $1\leq i\leq n$ and $1\leq j\leq m$. 

The morphism $\varphi$ is a morphism of $\B(E)$-comodule algebras, is the identity on the coinvariants elements $k$ of $Z$, and Proposition~\ref{prop_galois_object} ensures that the comodule algebra~$\B(E,F)$ has a bijective canonical map. Moreover, $Z$ is a faithfully flat Galois extension of $k$. Then by \cite[Remark 3.11]{S3} the morphism $\varphi$ is an isomorphism, and~$Z$ and~$\B (E,F)$ are isomorphic~$\B(E)$-Galois objects.

It remains to prove that the size $m$ of $F$ $\geq 2$. First assume that~$m=1$. Then~$W =\omega _Z (V_E) \cong k$. By \cite{S1}, \cite{S2}, there is an Hopf algebra~$K$ such that~$Z$ is a~$\B(E)$-$K$-bi-Galois object. Since there exists an inverse~$Z^{-1}$ of~$Z$ for the groupoid structure of bi-Galois objects, we have 
\begin{equation} 
V_E \cong V_E \Box _{\B(E)} Z \Box _{K} Z^{-1} \cong k \Box _K Z^{-1}.  \non \end{equation} 
Since $Z^{-1}$ is a Galois object, the image~$k \Box _K Z^{-1}$ of the trivial comodule of dimension one is the algebra~$k\cong (Z^{-1})^{co H}$ of coinvariants. 
Then the size~$m$ of~$F$ is equal to one only if the size~$n$ of$~E$ is one. The same argument proves that $m$ cannot be zero.

\end{proof}

We now turn to the classification of the Galois objects~$\B(E,F)$. The following lemma, implicit in~\cite{B}, will be useful. 

\begin{lemma} \lbl{lemma_basis} Let $k$ be a PID, let $n,m\geq 2$ be integers, and $E \in GL_n(k)$ and $F\in GL_m(k)$ be invertible matrices. 
Assume that~$\B(E,F)$ is a~$\B(E)$-$\B(F)$-bi-Galois object and let~$\varphi : \Comod $-$\B(E) \to \Comod $-$\B(F)$ be the associated monoidal equivalence. Let~$V_E$ and~$V_F$ be the respective fundamental comodules of~$\B(E)$ and~$\B(F)$. Then 
\begin{equation} \varphi (V_E) \cong V_F \non  \end{equation} 

% and let~$V$ be a finite free~$k$-module of rank~$n$ endowed with the~$\B(E)$-comodule structure defined by~(\ref{def_comod_V}). Let~$(z_{ij})_{i=1,\ldots 
%,n; j=1,\ldots ,m_1}$ be the canonical of generators of~$\B(E,F)$ (see~(\ref{def_BEF})). 
%Then~$(\sum _{k=1}^n v_k \otimes z_{kj})_{j=1,\ldots ,m}$ is a basis of~$V \Box _{\B(E)} \B(E,F)$.
\end{lemma}
\begin{proof}
Let~$w_1,\ldots ,w_m$ be the canonical basis of~$V_F$. Then we have a~$\B(F)$-colinear map~$\theta _F : V_F \to \varphi (V_E)$ defined by 

\begin{equation}  \theta _F(w_j)= \sum _{i=1}^n v_i \otimes z_{ij}.\non 
\end{equation} 
Similarly, we have a~$\B(E)$-colinear morphism~$\theta _E : V_E \to \varphi ^{-1} (V_F)$ defined by 

\begin{equation} \theta _E(v_i) = \sum _{j=1}^m w_j \otimes t_{ji} ,\non  \end{equation} where the~$t_{ji}$'s are the generators of~$\B(F,E)$. It is easy to see that $\varphi(\theta_E) \circ \theta_F$ is 
the canonical isomorphism $V_F \to \varphi (\varphi ^{-1} (V_F))$. We deduce that~$\theta _F$ and~$\theta _E$ are monomorphisms and then that~$\theta _F$ is an isomorphism. 
\end{proof}

As an immediate consequence of Lemma~\ref{lemma_basis}, we have the following necessary condition for~$\B(E)$-Galois objects to be cleft. 
\begin{cor}\lbl{cor_clive}
Let~$k$ be a PID and $n,m\geq 2$ be integers, $ E \in GL_n(k)$, $F\in GL_m(k)$ and $\B(E,F)$ be a cleft~$\B (E)$-Galois object. Then $m=n$. 
\end{cor}
\begin{proof}
If $\B(E,F)$ is a cleft Galois object, the associated fibre functor is isomorphic as a functor to the forgetful functor and in particular preserves the rank of finite free modules. 
\end{proof}

Let us now state our classification result for the extensions~$\B(E,F)$. 

\begin{theorem}\lbl{classif_F1F2}Let $k$ be a PID,~$n,m_1,m_2  $ be integers~$\geq 2$ and $E\in  GL _n (k) ,F_1 \in GL_{ m_1}(k) , F_2 \in GL_{m_2}(k)$ be invertible matrices such that the algebras~$\B(E,F_1)$ and~$\B(E,F_2)$ 
are~$k$-faithfully flat. Then the~$\B(E)$-Galois objects~$\B(E,F_1)$ and $\B(E,F_2)$ are isomorphic if and only if~$m_1=m_2$ and there 
exists an invertible matrix~$P\in GL_{m_1}(k)$ such that~$F_1 =PF_2 P^t$. 
\end{theorem}

Note that, by \cite{R} the bilinear forms associated to~$F_1$ and~$F_2$ are equivalent if and only if the asymmetries of~$F_1$ and~$F_2$ are similar. 

\begin{proof}As in the proof of \cite[Proposition 2.3]{B}, one shows that if $P\in GL_m(k)$, the~$B(E)$-comodule algebras $\B (E,F)$ and $\B(E,PFP^t)$ are isomorphic. 
 
Conversely assume that $\B(E,F_1)$ and $\B(E,F_2)$ are~$k$-faithfully flat: then Proposition \ref{prop_galois_object} ensures 
that~$\B(E,F_1)$ and~$ \B(E,F_2)$ are Galois objects. Let $V_E$ be the fundamental~$\B(E)$-comodule and let~$\beta _E : V_E \otimes V_E \to k $ be the linear map defined by~$E$. Let~$\omega _1 =  - \BB \B(E,F_1)$ and $\omega _2 = - \BB \B(E,F_2)$ be the fibre functors associated to~$\B(E,F_1)$ and~$\B(E,F_2)$. 

By Lemma \ref{lemma_basis}, the vector space $\omega _1 (V_E)$ has a basis~$(w^1_1, \ldots ,w^1_{m_1})$ and~$\omega _2 (V_E)$ has a basis~$(w^2_1,\ldots ,w^2_{m_2})$. The comodule algebra isomorphism $\varphi :  \B(E,F_1) \to \B(E,F_2)$ induces an isomorphism~$\id \otimes \varphi : \omega _1 (V_E) \to \omega _2 (V_E)$. Then in particular the rank of these two free~$k$-modules is the same, that is~$m_1=m_2=m$. Let~$P \in GL_m (k)$ be the matrix of~$\id \otimes \varphi$  in the bases~$(w^1_1, \ldots ,w^1_{m_1})$ and~$(w^2_1,\ldots ,w^2_{m_2})$. 

The matrices of the bilinear maps~$\omega_1 (\beta _E)$ and~$\omega_2(\beta _E)$, in the bases~$(w^1_1, \ldots ,w^1_{m})$ and~$(w^2_1,\ldots ,w^2_{m})$, are~$F_1$ and~$F_2$ respectively. Moreover, the isomorphism~$\varphi$ gives the relation 
\begin{equation} \omega _1(\beta _E) = \omega _2(\beta _E) \circ ( ( \id \otimes  \varphi) \otimes (\id \otimes \varphi )). \non  \end{equation} 
That is for any $i,j =1,\ldots ,m$ 
\begin{equation} \begin{array}{rcl}
\omega_1 (\beta _E) ( w^1_i \otimes  w^1_j ) &=& \omega _2 (\beta _E) \left( (\sum _{k=1}^m P _{ik} w^2_k )\otimes  (\sum _{l=1}^m P _{jl} w^2_j ) \right) \\ 
(F_1)_{ij} &=& \sum _{k,l =1 }^m P_{ik} P_{jl} (F_2)_{kl}, \\ 
\end{array} \non 
\end{equation} or in matrix form $F_1 = P F_2 P^t$. 
\end{proof}

\begin{rmk}{\rm 
As an application of Theorem \ref{classif_F1F2}, let us consider the case where the matrix $F$ is symmetric. Let $k$ be a PID, let $n,m,p\geq 2$ be integers, and $E \in GL_n(k)$, $F\in GL_m(k)$ and $G \in GL _p(k)$ be invertible matrices. Assume that $F$ is symmetric and $\B(E,F)$ is a Galois object. Then the Galois objects $\B(E,F)$ and $\B(E,G)$ are isomorphic if and only if $G$ is symmetric of size $p=m$. }
\end{rmk}

We now consider the case when $k$ is a field. For any integer~$n\geq 2$, and any invertible matrix~$E\in GL _n (k)$ we define
\begin{equation} 
X_0 (E) =\{ F \in GL _m (k), m \geq 2, \Tr (F^{-1 }F^t) = \Tr (E^{-1}E^t) \}. \non 
\end{equation} 
Consider the equivalence relation $\sim$ defined by $F_1\sim F_2$ if and only if there exists~$P \in GL (k)$ such that~$F_1= P F_2 P^t$ and put~$X(E) = X_0 (E)/\sim $. 

\begin{cor}Assume that~$k$ is a field. Then for any~ $n\geq 2$ and~$E \in GL _n(k)$, there is a bijection $\psi :  X (E)  \to  \Gal (\B(E))$ sending~$F$ onto~$[\B(E,F)]$. 
\end{cor}
\begin{proof}
Propositions 3.2, 3.3 and 3.4 in \cite{B} ensure that we have indeed this map~$\psi$. Moreover, $\psi$ is surjective by Theorem \ref{theorem_classif_iso} and injective by Theorem \ref{classif_F1F2}. 
\end{proof} 

We also have the following result. 

\begin{cor}Assume that~$k$ is an algebraically closed field of characteristic zero. For any~$n\geq 2$ and~$E \in GL _n(k)$, the group of~$\B(E)$-$\B(E)$-bi-Galois objects is trivial. 
\end{cor}

\begin{proof}
Let~$Z$ be a~$\B(E)$-$\B(E)$-bi-Galois object. By Theorem \ref{theorem_classif_iso}, there exist~$m \geq 2$ and~$F \in GL _m (k)$ such that~$Z$ is isomorphic to~$\B(E,F)$ as a~$\B(E)$-Galois object. Bichon \cite[Propositions 3.3, 3.4]{B} has proved that~$\B(E,F)$ is a~$\B(E)$-$\B(F)$-bi-Galois object. Then by \cite[Theorem 3.5]{S1} the Hopf algebras~$\B(E)$ and~$\B(F)$ are isomorphic that is, by \cite[Theorem 5.3]{B}, there exists~$P\in GL (k)$ such that~$F= P^t EP$. The matrix~$P$ enables us to construct an isomorphism of left~$\B(E)$-Galois objects $Z\cong \B(E,F) \cong \B(E)$.
 Now since~$Z$ is a~$\B(E)$-bi-Galois object, we know from \cite[Theorem 3.5]{S1} that there exists~$f \in \mathrm{Aut} (\B(E))$ such that $Z\cong \B(E) ^f$ as $\B(E)$-bi-Galois objects. Such a bi-Galois object is trivial if and only if~$f$ is coinner. Since by \cite[Theorem 5.3]{B} any Hopf automorphism of~$\B(E)$ is coinner, we are done. 
\end{proof}

The lazy cohomology group of a Hopf algebra was introduced in \cite{BC}, where it was realized as a subgroup of the group of bi-Galois objects. Therefore, we have the following. 
\begin{cor}Assume that~$k$ is an algebraically closed field of characteristic zero. The lazy cohomology group of~$\B(E)$ is trivial for any~$E\in GL _n(k)$. 
\end{cor} 

\section{Galois objects up to homotopy}

In this section we study the homotopy theory of~$\B(E)$-Galois objects. We assume that $k$ is an algebraically closed field. For technical reasons we only consider~$\Oq$-Galois objects (recall that~$\Oq =\B(E_q)$). Since for any~$E\in GL_n (k)$ there exists a~$\B(E)$-$\B(E_q)$-bi-Galois object, Proposition~\ref{prop_map} ensures that~$\h (\B(E)) \cong \h (\B(E_q) )$ and then there is no loss of generality. 

%then, \cite[Section 4]{B}, the Galois objects $\B(E_q,F)$ are nonzero algebra provided~$\mathrm{Tr}(F^{-1\;t}F) =-q-q^{-1}$. 
%In the case where the Hopf algebra $\B(E)$ is the quantum group $\Oq$ the previous classification Theorem \ref{theorem_classif_iso} could be 
%summerized in the following. 

We begin, using Lemma \ref{lemma_basis}, by giving a necessary condition for two~$\B(E_q)$-Galois extensions to be homotopically equivalent.

\begin{prop}\lbl{prop_neces}
Let~$m_0,m_1 \geq 2 $ be integers, let $F_0 \in GL_{m_0} (k) ,$ $F_1 \in GL_{m_1}(k)$ and assume that~$\B(E_q,F_0)$ and~$\B(E_q,F_1)$ are~$\B(E_q)$-Galois objects. If~$\B(E_q,F_0)$ and~$\B(E_q,F_1)$ are homotopically equivalent, then the matrices~$F_0$ and~$F_1$ have the same size $m_0=m_1$.
\end{prop}
\begin{proof}

%Let us consider a free finetely generated $k$-module~$V$ of rank~$n$ with the~$\B(E)$-comodule structure defined by (\ref{def_comod_V}). Then the 
%cotensor product define a map $\nu _V : \Gal _k (\B(E)) \to \N$ by 
%\begin{equation} \nu _V( Z)= \rk (V \Box _{\B(E)} Z ), \end{equation} for any Galois object~$Z$, where $\rk$ denotes the rank as~$k$-module, which 
%induces a map~$\tilde{\psi} : \h _k (\B(E)) \to \N$. 
%By Lemma \ref{lemma_basis}, a basis of~$V\Box _{\B(E)} \B(E,F)$ is~$(\sum _{k=1}^n v_k \otimes z_{kj})_{j=1\ldots m}$, where~$(v_k)_{k=1\ldots n}$ is a 
%basis of~$V$ and~$(z_{ij})_{i=1\ldots n; j=1\ldots m}$ is a set of generators of~$\B(E,F)$. Then the rank~$\rk(V\Box _{\B(E)} \B(E,F))$ is equal to the 
%dimension of the matrix~$F$. 

%let us now consider two~$\B(E)$-Galois objects~$\B(E,F_0)$ and~$\B(E,F_1)$ with homotopy~$\B(E,F)$. Then the rank of~$\B(E,F_0)$ and~$\B(E,F_1)$ are 
%respectively the dimensions of the matrices~$F_0$ and~$F_1$ and the rank, as~$k[t]$-module, of~$V\Box _{\B(E)} \B(E,F)$ is the dimension of the matrix~
%$F\in GL (k[t])$ which doesn't depend on~$t$. When~$t=0$ or~$t=1$, the dimension of~$F$ is the one of~$F_0$ and~$F_1$ respectively, then thehave the 
%same size~$m_0=m_1$. 
%By Lemma \ref{lemma_basis}, $V_E \Box _{\B(E)} \B(E,F)$ is a finite free~$k$-module of rank equal to the size of the matrix~$F$. 

Let us consider two~$\B(E_q)$-Galois objects~$\B(E_q,F_0)$ and~$\B(E_q,F_1)$ with homotopy~$\B_{k[t]}(E_q,F_t)$ (by Theorem \ref{theorem_classif_iso}, any~$\B_{k[t]}(E_q)$-Galois object is of this form for some~$F_t\in GL_m(k[t])$). Then~$V_E \Box _{\B_{k[t]}(E_q)} \B_{k[t]}(E_q,F_t)$ is a finite free~$k[t]$-module of rank equal to the size of the matrix~$F_t$, which does not depend on~$t$. The evaluation at~$t=0,1$ gives~$m_0=m_1$.
\end{proof}

%Note that we have a more general result. 
%\begin{prop}
%Let~$k$ be a principal commutative ring, $H$ be an Hopf algebra and~$V$ a finite free~$H$-comodule. The map~$\varphi _V : \h _k (H) \to \N $ given by 
%\begin{equation} \varphi _V (Z)= \rk (V \Box _H Z), \end{equation} for any~$H$-Galois object~$Z$, is well defined. 
%\end{prop}

%\begin{proof}
%Let~$Z_0,Z_1$ be homotopycally equivalent $H$-galois objects and~$\omega _0, \omega _1$ be the associated fibre functors. Let define~$[i] : \mod (k[t]) 
%\to \mod (k) $ by~$[i] (V) = V\otimes _{k[t]} k $, where the~$k[t]$-module structure on~$k$ is given by~$t^n v=i^n v$, for any $i\in \N ,v\in V$. Note that~$[i]$ 
%maps finite free~$k[t]$-module of rank~$n$ on finite free~$k$-module of rank~$n$. As the Galois extensions are homotopically equivalent, there exists a 
%fibre functor~$\omega : \Comod (H) \to \mod (k[t])$ such that~$\omega _i = [i] \circ \omega$. Then the ranks~$\rk (V \Box _H Z_0)$ and~$\rk (V \Box _H Z_1)
%$ are equal to the rank of the~$k[t]$-module~$V\Box _H Z$.   

%\end{proof}

Let us state a sufficient condition for two~$\B(E)$-Galois objects to be homotopically equivalent. 

\begin{theorem}\lbl{theorem_classif_hom}Let~$k$ be an algebraically closed field, $m_0,m_1\geq 2$ be integers and~$F_0,F_1$ be invertible matrices of size~$m_0,m_1$ such that~$\Tr (F_i^{-1}F^t_i)=-q -q^{-1}$ for~$i=0,1$. 

If~$m_0 =m_1$ and if~$F_0^{-1}F_0^t$ and~$F_1^{-1} F_1^t$ have the same characteristic polynomial, then the two~$\Oq$-Galois objects~$\B (E_q,F_0)$ and~$\B(E_q,F_1)$ are homotopically equivalent.
\end{theorem}

The rest of the section is devoted to the proof of the theorem. To this end, we construct a homotopy between the Galois objects, that is an~$\Oq [t]$-Galois object over the polynomial ring $k[t]$. First, let us begin with some terminology. We will say that a matrix~$F \in GL _m (k)$ (here~$k$ is an arbitrary commutative ring) is {\bf manageable} if $F^{-1}_{mm}=0$ and if the rightmost nonzero coefficient~$F^{-1}_{mv}$ in the bottom row is an invertible element of~$k$. In the case of a manageable matrix, the proof of \cite[Proposition 3.4]{B} still works and we obtain: 
\begin{prop}
Assume that~$k$ is a commutative ring and let~$F\in GL _m(k)$ be a manageable matrix such that $\Tr (F^{-1}F^t)= -q-q^{-1}$. Then~$\B(E_q,F)$ is a free $k$-module.
\end{prop}

%Let us give explicit forms for the matrices involved in this construction and in particular the matrix $F_t$ which gives the homotopy.
% and then we will prove that $\B(E_q,F_t)$ is $\Oq$-Galois extension of $k[t]$ over $k[t]$ and then by \cite[Proposition 1.2]{KS}, it is faithfully flat. 
%By work of Riehm \cite{R}, the eigenvalues of an asymmetry can be~$1$ in a Jordan block of odd dimension, $-1$ in a Jordan block of even dimension, 
%or~$q\in k^*$ in a Jordan block of dimension $n$ and then $q^{-1}$ is also an eigenvalue of Jordan block of dimension $n$. 

The problem for constructing an homotopy is the following one. 

\begin{Prob}
Let~$F_0,F_1 \in GL_m(k)$ be manageable matrices such that $\Tr (F_0^{-1}F_0^t) =\Tr (F_1^{-1}F_1^t)$.  
Find a matrix~$F(t) \in GL _m (k[t])$ such that 
\begin{enumerate}
\item $F(0)=F_0, \; F(1)=F_1$. 

\item $\Tr (F(t)^{-1}F(t)^t) =\Tr (F_0^{-1}F_0^t) =\Tr (F_1^{-1}F_1^t)$. 

\item $F(t)$ is manageable. 
\end{enumerate}
\end{Prob}

Now assume that~$F_0$ and~$F_1$ have diagonal block decompositions with the same size: 
\begin{equation} 
F_0 =\left( \begin{array}{c|c}
(F_0) _{11} & 0 \\ 
\hline 
0 & (F_0)_{22} \\ \end{array} \right) , 
\quad F_1 =\left( \begin{array}{c|c}
(F_1) _{11} & 0 \\ 
\hline 
0 & (F_1)_{22} \\ \end{array} \right),\non 
\end{equation} that 
\begin{equation}
 \Tr ((F_0)_{11}^{-1}(F_0)_{11}^t) =\Tr ((F_1)_{11}^{-1}(F_1)_{11}^t) \non 
\end{equation} and
\begin{equation} 
\Tr ((F_0)_{22}^{-1}(F_0)_{22}^t) =\Tr ((F_1)_{22}^{-1}(F_1)^t_{22})\non   \end{equation} 
and finally that each block is manageable. Then clearly Problem~{\bf {\sf (P)}} reduces to the same problem for each block. This simple remark, combined with Riehm's work \cite{R} on the structure of bilinear forms, will reduce our problem to the case of some ``elementary" matrices. 

We will use freely the following results of \cite{R}. For any nondegenerate bilinear map~$\beta : V \times V \to k$ given by an invertible matrix~$F$, and for any eigenvalue~$p \neq \pm 1$ of its asymmetry~$\sigma$, $p^{-1}$ is also an eigenvalue of~$\sigma$ and the two characteristic spaces~$C_p$ and~$C_{p^{-1}}$ associated to~$p$ and~$p^{-1}$ are isotropic (for the bilinear form~$\beta$). The vector space~$V$ is the orthogonal sum of the subspaces~$C_1,C_{-1}$ and $C_p\oplus C_{p^{-1}},$ where~$p$ runs over all eigenvalues of~$\sigma$ different from~$\pm 1$. Then there exists a basis of~$V$ such that the matrix of~$\sigma$ is a block matrix made of Jordan blocks of odd dimension with eigenvalue~$1$, Jordan blocks of even dimension with eigenvalue~$-1$ and pairs of blocks of eigenvalues~$p,p^{-1}$ and of the same dimension. 

Assume that the asymmetries~$\sigma _0$ and ~$\sigma _1$ associated to~$F_0$ and~$F_1$ have the same characteristic polynomial and that~$\sigma _1$ is diagonal. Then by \cite{R}, Problem~{\textbf{\sf (P)}} reduces to three cases. 

\begin{itemize}
\item[A.] $\sigma _0$ is a Jordan block of even dimension~$d$ with eigenvalue~$-1$ (and~$\sigma _1 = - \mathrm{I}_d$),

\item[B.] $\sigma _0$ is a Jordan block of odd dimension~$d$ with eigenvalue~$1$ (and~$\sigma _1 =\mathrm{I}_d$),

\item[C.] $\sigma _0$ is a diagonal block matrix made of two Jordan blocks of eigenvalues~$p,p^{-1}$ and of the same size~$d$ (and~$\sigma _1$ is diagonal with $d$ diagonal coefficients equal to~$p$ and~$d$ equal to~$p^{-1})$. 
\end{itemize}

%Let~$\sigma \in GL_n (k)$ be the asymmetry of the bilinear map~$\beta : V_E \otimes V_E \to k$ and let 
%\begin{equation} V_p =\{ v\in V | p^r v=0 , \text{ for } r>>0 \} \end{equation} 
%for each irreductible divisor~$p$ of the caracteristic polynom of~$\sigma$. Then (\cite{R})~$V_E$ is the orthogonal sum of the subspaces 
%\begin{equation} V_p \oplus V_{p^{-1}} \quad (p\neq \pm 1), \quad V_1,\quad V_{-1} \end{equation}
%and there exists a basis such that the matrix of~$\sigma$ is a block matrix made with Jordan blocks of eigenvalue~$1$ and odd dimension, Jordan 
%blocks of eigenvalue~$-1$ and even dimension and pair of blocks of eigenvalues~$q,q^{-1}$ and of the same dimension. 

Let us now look at the possible forms of a matrix~$F$ such that~$\sigma _0=F^{-1}F^t $ for each of these three cases. 

\begin{lemma}\lbl{lemme_form_F}
A) If $\sigma _0$ is a Jordan block of even dimension with eigenvalue equal to~$-1$ and if there exists an invertible matrix $F$ such that $F^{-1}F^t = \sigma _0$, then $F$ has the lower anti-triangular form 
\begin{equation} \lbl{form_F_-}F= \left( \ba{ccccc} 0 &\cdots &\cdots & 0 &F_{1n} \\ \vdots & & \nedots & -F_{1n} &* \\\vdots & \nedots & \nedots &* &* \\0 & F_{1n} &* &* &* \\ -F_{1n} &*&*&*&* \\ \end{array} \right). \end{equation}

B) If $\sigma _0$ is a Jordan block of odd dimension with eigenvalue equal to $1$ and there exists an invertible matrix~$F$ such that~$F^{-1}F^t = \sigma _0$, then~$F$ has the lower anti-triangular form  
\begin{equation} \lbl{form_F_+}F= \left( \ba{ccccc} 0 &\cdots&\cdots & 0 &F_{1n} \\
 \vdots &&\nedots & - F_{1n} &* \\
\vdots &\nedots &\nedots &*&*\\
 0  & - F_{1n} &* &*&* \\  
F_{1n} &* &* &*&*  \\ \end{array} \right). \end{equation}

C) If $\sigma _0$ is made of two Jordan blocks of eigenvalues $p$ and $p^{-1}$ and of size~$n$, then the invertible matrix $F$ defined by 
\begin{equation}\lbl{def_F_pp}
\left( \begin{array}{cc} 
0 & \mathrm{I_n} \\
\mathrm{J_p} & 0 \\ 
\end{array}\right), 
\end{equation} where~$\mathrm{I_n}$ is the identity of size $n$, $\mathrm{J_p}$ is a Jordan block of size~$n$ and eigenvalue~$p$ and~$0$ is the zero matrix, has an asymmetry similar to~$\sigma _0$. 
\end{lemma}

\begin{proof}We say that the elements~$a_{i,n+1-i}$, for $1\leq i\leq n$ of a matrix~$A\in \mathrm{M} _n(k)$ lies on the anti-diagonal and we use obvious notion of lower and upper anti-triangular matrices. 

A) Assume that~$F$ is a matrix such that~$F^{-1}F^t=\sigma _0$, that is~$F^t = F \sigma _0$ or 
\be \lbl{equa_F_s-} \left\{ \begin{array}{ccll} F_{i1} &=&-F_{1i}& \forall i =1,\ldots ,n \\
F_{ji}&=&F_{i,j-1}-F_{ij} & \forall i=1,\ldots ,n; j=2,\ldots , n \\ \end{array}\right. \end{equation}
Let us consider the first row. The equation $F_{11}=-F_{11}$ implies~$F_{11}=0$. Then, for any $k=2,\ldots ,n$, we have from (\ref{equa_F_s-}) the equations $F_{1k}=-F_{k1}$ and~$F_{k1}= F_{1, k-1}-F_{1k}$ and then $F_{1,k-1}=0$. %, as $F_{1,k-1}=-F_{k-1,1}$. 
Then the first row and the first column are equal to zero except the last terms $F_{1n}$ and $F_{n1}$. 

For the second row and column, we have from the previous computation $F_{12}=F_{21}=0$ then $F_{22}=F_{21} -F_{22}=0$. For any $k=3,\ldots ,n-1$ we have 
\begin{equation} \left\{ \begin{array}{lll} F_{2k}&=&F_{k1}-F_{k2} \\ F_{k2}&=&F_{2,k-1}-F_{2k}. \\\end{array} \right.\non  \end{equation} 
Then $F_{2,k-1}=0$ and, since $F_{2,k-1}=F_{k-1,1}-F_{k-1,2}$, we also have $F_{k-1,2}=0$. Then for all $k \leq n-2$ the entries $F_{2,k}$ and $F_{k,2}$ are equal to zero. In the same way, any coefficient lying above the anti-diagonal is equal to zero. 

The coefficient $F_{i,n+1-i}$ on the anti-diagonal satisfies the relation $F_{i,n+1-i}= F_{n+1-i,i-1} - F_{n+1-i,i}= - F_{n+1-i,i}$. We also have
\begin{equation} \lbl{rel_diag1} \left\{ \begin{array}{l} F_{i,n+2-i}= F_{n+2-i,i-1}-F_{n+2-i,i} \\
F_{n+2-i,i} = F_{i,n+1-i}-F_{i,n+2-i},\\ \end{array} \right. \end{equation} 
then 
\begin{equation} F_{n+2-i,i}= F_{i,n+1-i}-F_{n+2-i,i-1}+ F_{n+2-i,i}   \end{equation} that is 
\begin{equation} \lbl{rel_diagder}F_{i,n+1-i}=F_{n+1-(i-1),(i-1)}. \end{equation}

Then the determinant of $F$ is $(F_{1n})^n$ and the matrix has the wanted form. 
%\begin{equation} F= \left( \ba{ccccc} 0 &\cdots &\cdots & 0 &F_{1n} \\ \vdots & & \nedots & -F_{1n} &* \\\vdots & \nedots & \nedots &* &* \\0 & 
%F_{1n} &* &* &* \\ -F_{1n} &*&*&*&* \\ \end{array} \right) \end{equation}

B) Let us now consider the case where $\sigma _0$ is a Jordan block of odd dimension and eigenvalue~$1$ and~$F$ is a matrix such that~$F^{-1}F^t=\sigma _0$, that is~$F^t = F \sigma _0$ or 
\be \lbl{equa_F_s+} \left\{ \begin{array}{ccll} F_{1i} &=&F_{i1}& \forall i =1,\ldots ,n \\
F_{ji}&=&F_{i,j-1}+ F_{ij} & \forall i=1,\ldots ,n; j=2,\ldots , n \\ \end{array}\right. \end{equation}

Let us consider the first line. For any $k=2,\ldots ,n$ we have from (\ref{equa_F_s+}) the equations $F_{1k}=F_{k1}$ and $F_{k1}= F_{1, k-1}+F_{1k}$ and then $F_{1,k-1}=0$, since $F_{1,k-1}=F_{k-1,1}$, the first line and the first column are equal to zero except the last terms $F_{1n}=F_{n1}$. 
In the same way as for the previous case, we see that all the coefficients lying above the anti-diagonal must be zero. Moreover, in the same way as for (\ref{rel_diag1}) - (\ref{rel_diagder}), the anti-diagonal coefficient in position $(i,n-i+1)$ is $(-1)^{i+1}F_{1n}$; the determinant is~$(F_{1n})^n$ and~$F$ has the wanted form.
%\begin{equation} F= \left( \ba{ccccc} 0 &\cdots &\cdots & 0 &F_{1n} \\ \vdots & & \nedots & -F_{1n} &* \\\vdots & \nedots & \nedots &* &* \\0 
%&- F_{1n} &* &* &* \\ F_{1n} &*&*&*&* \\ \end{array} \right). \end{equation}

C) Assume that $\sigma _0$ is made of two Jordan blocks of eigenvalues~$p$ and~$p^{-1}$ and of size~$n$. We define~$F$ by the relation~(\ref{def_F_pp}). Its asymmetry is the matrix 
\begin{equation}\left( \begin{array}{cc} 
\mathrm{J_p}^{-1} &0 \\
0 &  \mathrm{J_p} ^t \\ 
\end{array}\right) \non 
\end{equation} which is similar to~$\sigma _0$.
\end{proof}

\begin{proof}[Proof of Theorem \ref{theorem_classif_hom}]
Let us construct the matrix~$F(t)$ solution of the problem~{\bf \sf (P)}. 

%\begin{lemma}\lbl{construction_F_t} 
%Let $F \in GL _n(k)$ be an invertible matrix and $\sigma$ be its asymmetry. There exists a matrix $F_t \in GL _n (k[t])$ such that the matrix $F_t$ 
%evaluated at $t=1$ is equal to $F$ and that its evaluation at $t=0$ is an invertible matrix $F_0$ with diagonal asymmetry $\sigma '$ and such that the 
%list of eigenvalues (with multiplicity) of $\sigma '$ is the same as the one of $\sigma$. Moreover the matrix~$F(t)$ is manageable and satisfy 
%\end{lemma}
%\begin{proof}
%By block computation, it is enough to look at the Jordan blocks of $\sigma $ and then at the two following cases : $\sigma$ is a Jordan block with 
%eigenvalues equal to $\pm 1$ or $\sigma$ is made of two Jordan blocks with the same size and eigenvalues equal to $q$ and $q^{-1}$. 

{\it Cases A, B:} Let us consider a Jordan block $\sigma _0$ with eigenvalue $\pm 1$ and size more than two. 
By the previous lemma \ref{lemme_form_F}, the matrix~$F$ such that~$F^{-1}F^t = \sigma _0$ is an anti-triangular matrix of form (\ref{form_F_-}) or (\ref{form_F_+}). Consider the matrix $F(t) \in GL _n(k[t])$ defined by 
\begin{equation}
F(t)_{i,n+1-i}=F_{i,n+1-i}, \quad F(t)_{ij}=tF_{ij},\non  \end{equation}
for any $1\leq i,j \leq n$ such that $j\neq n+1-i$ (that is $F(t)$ is equal to $F$ on the anti-diagonal and to $tF$ on the other coefficients).
 
To compute $\Tr (F(t)^{-1}F(t)^t)$ we have to know the diagonal coefficients of the asymmetry of $F(t)$, which are equal to products of the anti-diagonal coefficients of~$F(t)^{-1}$ and~$F(t)^t$.
Remark that if a matrix $F(t)$ is lower anti-triangular, its inverse $F(t)^{-1}$ is upper anti-triangular, and their anti-diagonal coefficients are related  by 
\begin{equation}
1= (F(t)^{-1})_{i,n+1-i} (F(t))_{n+1-i,i},\non  \end{equation}
for any $i=1, \ldots ,n$. Since the anti-diagonal coefficients of~$F(t)$ do not depend on $t$, the ones of $F(t)^{-1}$ do not depend on~$t$ either and we have 
\begin{equation} 
\Tr(F(t)^{-1}F(t)^t) =\Tr(F_0^{-1}F_0^t) =\Tr(F_1^{-1}F_1^t) .\non 
\end{equation}

From the definition of $F(t)$, we have $F(0)=F_0$ and $F(1)$ is a block matrix with anti-diagonal blocks. Then the asymmetry of $F(1)$ is diagonal and then equal to~$\sigma _1$. 
%we have seen that the diagonal coefficient of the asymmetry doesn't depend on $t$ and the are the eigenvalues of $\sigma _0$. 

Since $F(t)^{-1}$ is upper anti-triangular and invertible, $F(t)$ is manageable. Finally, $F(t)$ is a solution of {\bf \sf (P)} in the cases A,B. 

{\it Case C:} Let us now consider the case of two Jordan blocks of size $n$ and eigenvalues $p$ and $p^{-1}$ and suppose that $F$ has the form~(\ref{def_F_pp}). 

Consider the matrix $F(t) \in GL_{2n} (k[t])$ defined by 
\begin{equation}\left( \begin{array}{cc} 
0 & \mathrm{I_n} \\
\mathrm{J_p(t)} & 0 \\ 
\end{array}\right), \non 
\end{equation}
where $\mathrm{J_p(t)}$ is the matrix with diagonal coefficients equal to $p$ and upper diagonal coefficients $(i,i+1)$ equal to $t$. 
%and note $G_t$ the matrix defined from $F$ by changing $g$ by~$g_t$. 
%An easy computation give that the asymmetry related to $G_t$ has the same diagonal coefficient as the one related to $G$ and then its trace is 
%constant. 
The inverse~$F(t)^{-1}$ is 
\begin{equation}\left( \begin{array}{cc} 
0 & (\mathrm{J_p(t)} )^{-1} \\
 \mathrm{I_n}& 0 \\ 
\end{array}\right) \non 
\end{equation} and then~$F(t)$ is manageable and the trace of its asymmetry is constant. Finally~$F(t)$ is a solution of { \bf \sf (P)} in the case C. 
\end{proof}

\begin{cor} 
All cleft $\Oq$-Galois objects are homotopically trivial. 
\end{cor}

\begin{proof}
Let $Z$ be a cleft Galois object of~$\Oq$. Then by Theorem \ref{theorem_classif_iso}, the Galois object~$Z$ is isomorphic to $\B(E_q,F)$ and by Corollary $\ref{cor_clive}$ the matrix $F$ is a $2\times 2$ matrix with trace equal to~$-q-q^{-1}$. 

If $q\neq 1$, there exists~$P\in GL (k)$ such that~$F=PE_q P^t$. Then the Galois object $\B(E_q,F)$ is isomorphic to the trivial object $\B(E_q)$. 

If $q= 1$, the two possible asymmetries are, up to similarity, a diagonal matrix $\sigma_1$ with eigenvalue $-1$ and multiplicity $2$ associated to a matrix $F_1$ or a Jordan block matrix $\sigma _2$ of size $2$ and eigenvalue $-1$ associated to a matrix $F_2$. The two associated Galois objects $\B(E_q,F_1)$ and $\B(E_q,F_2)$ are nonisomorphic as the asymmetries are nonsimilar, but they are homotopically equivalent by Theorem~\ref{theorem_classif_hom} as the asymmetries have the same characteristic polynomial. 
\end{proof}

\end{document}